\documentclass[12pt]{article}
\usepackage{amssymb}

\newtheorem{thm}{Theorem}[section]
\newtheorem{dfn}{Definition}[section]
\newtheorem{prop}{Proposition}[section]
\newtheorem{lem}{Lemma}[section]
\newtheorem{rem}{Remark}[section]
\newtheorem{cor}{Corollary}[section]

\def\B{{\cal{B}}}
\def\E{{\cal{E}}}

\def\C{{\mathbb{C}}}

\def\DD{{\bf{D}}}
\def\M{{\cal{M}}}

\def\Z{{\mathbb{Z}}}

\def\la{{\langle}}
\def\ra{{\rangle}}
\def\lla{{\langle \! \langle}}
\def\rra{{\rangle \! \rangle}}
\title{Braided differential structure on Weyl groups, quadratic algebras and elliptic functions 
\vspace {5mm} \\ 
\small\it{To the memory of Leonid Vaksman}} 
\date{}
\author{Anatol N. Kirillov and Toshiaki Maeno\footnote{Supported by Grant-in-Aid for Scientific Research.}}
\begin{document}
\maketitle
\begin{abstract}
We discuss a class of generalized divided difference operators which give 
rise to a representation of Nichols-Woronowicz algebras associated to Weyl 
groups. For the root system of type $A,$ we also study the condition for 
the deformations of the Fomin-Kirillov quadratic algebra, which 
is a quadratic lift of the Nichols-Woronowicz algebra, 
to admit a representation given by generalized divided difference operators. 
The relations satisfied by the mutually commuting elements called Dunkl
elements in the deformed Fomin-Kirillov algebra are determined. 
The Dunkl elements correspond to the truncated elliptic Dunkl operators 
via the representation given by the generalized divided difference 
operators. 
\end{abstract}

\section*{Introduction} 
The rational Dunkl operators, which were introduced in \cite{Du} for any 
finite Coxeter group, constitute a remarkable family of operators of differential-difference type. 
The Dunkl operators 
are defined to be the ones acting on the functions on the reflection 
representation $V$ of the corresponding Weyl group $W.$ For the root system of type $A_{n-1},$ 
the Dunkl operators $D_1,\ldots,D_n$ are defined by the formula 
\[ D_i:= \frac{\partial}{\partial x_i}+ 
\sum_{j \not = i}\frac{1-s_{ij}}{x_i-x_j} , \] 
where $s_{ij}$ is the transposition of $i$ and $j.$ 
They are $S_n$-invariant and mutually commute. The Dunkl opearotors play an 
important role in the 
representation theory and in the study of integrable systems. Here we would like 
to mention only a remarkable result, due to Dunkl, that the algebra 
generated by {\it truncated Dunkl operetors} is isomorphic to the coinvariant
algebra of the corresponding finite Coxeter group \cite{Du}, \cite{Ba}. A 
trigonometric generalization of Dunkl operators has been proposed by 
Cherednik \cite{Ch}, and an elliptic one by Buchstaber, Felder and Veselov 
\cite{BFV}. The basic requirement for such generalizations is that the 
operators to be constructed are bounded to pairwise commute. Another 
important property of rational Dunkl operators, namely, their 
$W$-invariance, may be broken for generalizations. 

For a crystallographic irreducible root system $R,$ Buchstaber, Felder and 
Veselov \cite{BFV} have determined the conditions on the functions 
$f_{\alpha}(z),$ $\alpha\in R,$ so that the operators 
\[ \nabla_{\xi}=\partial_{\xi}+\sum_{\alpha\in R_+}(\alpha,\xi)f_{\alpha}((\alpha,x))
s_{\alpha} \] 
satisfy the commutativity condition $[ \nabla_{\xi},\nabla_{\eta} ]=0$ for 
all $\xi,\eta \in V.$ Here, we denote by $R_+$ the set of positive roots and 
by $s_{\alpha}$ the reflection corresponding to a root $\alpha.$ 
Under the assumption of the $W$-invariance of $\nabla_{\xi},$ they proved that 
the solutions of the functional equation for $f_{\alpha}$ must be rational 
unless $R$ is of type $B_2.$ Without the assumption of the $W$-invariance, 
some elliptic solutions given by Kronecker's 
$\sigma$-function may appear. If $R$ is of type $A_n,$ such functions exhaust 
the general solution. 

The present paper contains two main results. The first one is concerned about 
the existence of a representation given by the generalized divided difference 
operators for the (certain extension of) Nichols-Woronowicz algebra $\B_W$ 
corresponding to a Weyl group $W.$ Our second main result describes relations 
among the Dunkl elements in the elliptic extension of the Fomin-Kirillov 
algebra introduced originally in \cite{K1}. 
In particular, we describe the relations among truncated elliptic 
Dunkl operators of type $A_{n-1}.$ 
By analogy with Dunkl's theorem 
mentioned above, one can consider the algebra generated by truncated elliptic 
Dunkl 
operators of type $A_{n-1}$ as an elliptic deformation of the cohomology ring 
of the flag variety $Fl_n.$ We also prove an elliptic analogue of the 
Pieri rule in the elliptic extension of Fomin-Kirillov algebra. These 
results can be 
considered as further generalizations of those obtained in \cite{FK}, 
\cite{Po}, since the latter correspond to certain degenerations of the 
elliptic case, see Section 4 for details. 

The Nichols-Woronowicz algebra $\B(M)$ is 
a braided analogue of the symmetric algebra, which is defined for a given 
braided vector space $M.$ Nichols \cite{Ni} studied graded bialgebras 
generated by the primitive elements of degree one. 
The braided Hopf 
algebra $\B(M)$ satisfying such a condition was called Nichols algebra 
by Andruskiewitsch and Schneider \cite{AS}. 
The algebra $\B(M)$ has 
been constructed also in the theory of the differential 
forms on quantum groups due to Woronowicz \cite{Wo}. 
Woronowicz constructed $\B(M)$ as a braided symmetric (or exterior) algebra 
based on the construction of his (anti-)symmetrizer. 
The Nichols-Woronowicz algebra provides a natural framework for the braided 
differential calculus, which was developed by Majid \cite{M1}. 

In this paper we are interested in the Nichols-Woronowicz algebra associated 
to a particular kind of braided vector space called Yetter-Drinfeld module. 
See \cite{Ba} for more details of general construction of $\B(M).$ 
In our case, we use a $\C$-vector space $M_W$ spanned by the symbols 
$[\alpha]=-[-\alpha],$ $\alpha \in R,$ with the braiding 
$\psi : M_W^{\otimes 2}\rightarrow M_W^{\otimes 2},$ 
$[\alpha]\otimes [\beta] \mapsto [s_{\alpha}(\beta)] \otimes [\alpha].$ 
The algebra $\B_W=\B(M_W)$ of our interest is defined to be the quotint of the tensor 
algebra of $M_W$ by the kernel of the braided symmetrizer. 

Milinski and Schneider \cite{MS} and Majid \cite{M2} have pointed out that 
the algebra $\B_W$ for 
$W=S_n$ is a quotient of the Fomin-Kirillov quadratic algebra 
$\E_n$ defined in \cite{FK}. 
The algebra $\B_{S_n}$ is conjectured to be isomorphic to $\E_n.$ 
Fomin and the first author introduced the 
algebra $\E_n$ to construct a model of the cohomology ring of the flag 
variety $Fl_n.$ In \cite{Ba}, Bazlov has reformulated their construction of 
the model of the cohomology ring in terms of the Nichols-Woronowicz algebra 
$\B_W,$ and generalized it to arbitrary finite Coxeter groups. The braided 
differential operators on the algebra $\B_W,$ which were used by Majid \cite{M2} 
for root system of type $A,$ play an essential role in Bazlov's 
construction. His construction also has an important implication on the 
representation of $\B_W,$ since the braided differential operators act on the 
coinvariant algebra of $W$ as the divided difference operators 
$\partial_{\alpha}=(1-s_{\alpha})/\alpha,$ $\alpha\in R.$ 

In Section 1, we discuss the conditions for the generalized divided 
difference operators 
\[ \DD_{\alpha}=f_{\alpha}((\alpha,\xi))+g_{\alpha}((\alpha,\xi))s_{\alpha} \]
to give rise a representation of $\B_W.$  
These conditions are interpreted as functional equations for $f_{\alpha}$ and $g_{\alpha}.$ 
We prove that the operators 
corresponding to the $W$-invariant solutions described in \cite{BFV} 
define a representation of $\B_W.$ Komori \cite{Ko} studied when the operators 
$\DD_{\alpha}$ satisfy the Yang-Baxter equation. Since the generators $[\alpha]$ 
of the algebra $\B_W$ satisfy the Yang-Baxter equation, our operators also correspond to 
special part of the solutions found in \cite{Ko}.

In order to get a more general class of solutions like elliptic functions, 
we have to loose part of defining relations of $\B_W.$ 
In Section 2 we introduce a deformed version 
$\tilde{\E}_n(\psi_{ij})$ of the Fomin-Kirillov quadratic algebra, 
which is defined for a given family of meromorphic functions 
$\psi_{ij}(z),$ $1\leq i,j \leq n,$ 
$i\not=j.$ The algebra $\tilde{\E}_n(\psi_{ij})$ admits the representation by 
the operators $\DD_{\alpha}$ only when $\psi_{ij}(z)$ is given by the 
Weierstrass $\wp$-function or its degenerations.  In this case, the operator 
$\DD_{\alpha}$ exactly corresponds to the general solution 
for $A_{n-1}$-system obtained in \cite[Theorem 16]{BFV}. 

Our second main result is the study of relations among the Dunkl elements in 
the elliptic extension $\tilde{\E}_n(\psi_{ij})$ of the Fomin-Kirilov algebra. 
The Dunkl elements $\theta_1,\ldots,\theta_n \in \tilde{\E}_n(\psi_{ij})$ 
are mutually commuting elements defined by $\theta_i=\sum_{j\not=i}[ij].$ 
The images of the Dunkl elements, via the representation 
$[\alpha] \mapsto \DD_{\alpha},$ become the so-called truncated 
(or level zero) elliptic Dunkl operators, cf \cite{BFV}. 
It is well-known that the (truncated) rational or trigonometric 
Dunkl operators can be obtained as certain degenerations of the (truncated) 
elliptic Dunkl operators. The identities among the Dunkl elements in $\tilde{\E}_n(\psi_{ij})$ 
are also satisfied by the corresponding truncated elliptic Dunkl operators 
or their degenerations. 
In the context of Schubert calculus, the Dunkl elements describe the 
multiplication by the classes of standard line bundles in the cohomology ring 
of the flag variety. The formula of the elementary symmetric polynomials 
in the Dunkl elements in the Fomin-Kirillov algebra reflects the Pieri 
formula. In Section 3 we give a formula for the deformed elementary symmetric 
polynomial $E_k(\theta_i \; | \; i\in I)$ in the algebra 
$\tilde{\E}_n(\psi_{ij}).$ 

The algebra $\tilde{\E}_n(\psi_{ij})$ 
has degenerations to variants of the deformation of the Fomin-Kirillov 
algebra. In particular, the multiparameter deformation $\E^p_n$ studied in 
\cite{FK} and \cite{Po}, 
and the extended quadratic algebra $\tilde{\E}_n \langle R \rangle [t]$ 
defined in \cite{KM2} after the specialization $t=0$ can be regarded as 
degenerations of $\tilde{\E}_n(\psi_{ij}).$ 
In Section 4 we show that our algebra recovers the Pieri formulas in the 
corresponding degenerations. 

\section{Representation of Nichols-Woronowicz algebra} 
Let us consider the reflection representation $V$ of 
the Weyl group $W.$ Denote by $R\subset V$ the set of roots for the Weyl group $W.$ Fix $R_+$ 
the set of 
positive roots in $R.$ Let $\{ \alpha_1,\ldots , \alpha_r \} \subset R_+$ be the set 
of simple roots. 
The Weyl group $W$ naturally acts on the space $\M=\M(V_{\C})$ of meromorphic functions on 
$V_{\C}.$ 
We also denote by $\M_0$ the space of meromorphic functions on $\C.$ 

We discuss the generalized Calogero-Moser representaion of the Nichols-Woronowicz algebra $\B_W$ 
for the Weyl group $W.$ The Nichols-Woronowicz algebra $\B_W=\B(M_W)$ is associated to the 
Yetter-Drinfeld module $M_W$ generated by the symbols $[\alpha],$ $\alpha \in R.$ Define the operator 
$\DD_{\alpha},$ $\alpha \in R,$ acting on $\M$ by 
\[ \DD_{\alpha}=f_{\alpha}((\alpha,\xi))+g_{\alpha}((\alpha,\xi))s_{\alpha}, \; \; 
\xi \in V, \]
where $s_{\alpha}$ is the reflection with respect to $\alpha,$ and $f_{\alpha},g_{\alpha}\in \M_0.$ 
We assume that $f_{-\alpha}(z)=-f_{\alpha}(-z)$ and $g_{-\alpha}(z)=-g_{\alpha}(-z)$ so that 
$\DD_{-\alpha}=-\DD_{\alpha}.$ \begin{lem} The divided difference operator 
$\partial_{\alpha}=(1-s_{\alpha})/(\alpha,\xi)$
gives a well-defined representation of $\B_W$ on $P.$ 
\end{lem} 
{\it Proof.} From the construction of 
the model of the coinvariant algebra $P_W$ in \cite{Ba}, 
we can see that the natural action of the braided differential operator 
$\overleftarrow{D}_{\alpha}$ on $P_W$ coincides with the divided difference operator 
$\partial_{\alpha}.$ Since $P=P^W \otimes P_W,$ we can extend $P^W$-linearly the action of 
$\B_W$ on $P_W$ to that on $P.$ 
\begin{lem} If $[\alpha] \mapsto \DD_{\alpha}$ defines the 
representation of $\B_W,$ then 
$f_{\alpha}$ must be an odd function, and \[ g_{\alpha}(z)=f_{\alpha}(z) \phi_{\alpha}(z), \] 
where $\phi_{\alpha}(z)\phi_{\alpha}(-z)=1.$
\end{lem}
{\it Proof.} The condition $\DD_{\alpha}^2=0$ is equivalent to the equations 
\[ f_{\alpha}(z)^2+g_{\alpha}(-z)g_{\alpha}(z)=0 \] and \[ g_{\alpha}(z)\cdot 
(f_{\alpha}(z)+f_{\alpha}(-z))=0 . \] The second equation shows that $f_{\alpha}$ is odd. 
Define the function $\phi_{\alpha}(z)$ by \[ \phi_{\alpha}(z)= \frac{g_{\alpha}(z)}{f_{\alpha}(z)} . \] 
Then the first equation can be written as \[ \phi_{\alpha}(z)\phi_{\alpha}(-z)=1. \] 

We take the standard realization of the root systems of type $A_n$ and $B_n$ as follows: 
\[ R(A_n)= \{ ij=\epsilon_i-\epsilon_j \; | \; 1\leq i,j \leq n, i\not= j \} , \] 
\[ R(B_n)= \{ ij=\epsilon_i-\epsilon_j, \overline{ij}=\epsilon_i+\epsilon_j, i=\epsilon_i | \; 
1\leq i,j \leq n, i\not= j \} , \] 
where $(\epsilon_1,\ldots ,\epsilon_n)$ is an orthonormal basis of $V.$ 
\begin{prop} Suppose that $R$ is not of type $A_1$ or $B_2.$ If the operators 
$\DD_{\alpha}$ give a representation of $\B_W,$ then 
$f_{\alpha}(z)= k_{\alpha}/z$ and $g_{\alpha}(z)=\pm k_{\alpha}e^{\lambda_{\alpha}z}/z,$
where $k_{\alpha}$ are $W$-invariant constants and the choice of the signature $\pm$ is 
independent of roots $\alpha.$ 
The constants $\lambda_{\alpha}$ are obtained as $\lambda_{\alpha}=\lambda(\alpha^{\vee})$
from an element $\lambda \in V^*.$
Conversely, the operators $\DD_{\alpha}$ corresponding to the above solutions 
give the representation of $\B_W.$ 
\end{prop} 
{\it Proof.} When $R$ is of type $A_2,$ we have the functional equations 
\begin{eqnarray}
f_{12}(x-y)f_{23}(y-z)+f_{23}(y-z)f_{31}(z-x)+f_{31}(z-x)f_{12}(x-y)=0, \\ 
g_{12}(x-y)g_{23}(x-z)+g_{23}(y-z)g_{31}(y-x)+g_{31}(z-x)g_{12}(z-y)=0. 
\end{eqnarray}
If $f_{12}$ is regular at the origin, then we have $f_{12}(0)=0$ since $f_{12}$ is odd. 
We have $f_{23}(x-z)f_{31}(z-x)=0$ 
by putting $x=y$ in the equation (1), and hence $f_{12},$ $f_{23}$ and $f_{13}$ must be constantly zero. 
So we may assume $f_{12},$ $f_{23}$ and $f_{13}$ have a pole at the origin. 
Now the equation (1) shows
\[ f_{31}(z-x)^{-1}+ f_{12}(x-y)^{-1} + f_{23}(y-z)^{-1}=0 . \] 
Therefore we have \[ f_{12}(x)=f_{23}(x)=f_{13}(x)= \frac{k}{x} \] 
for some constant $k.$ From Lemma 1.2, we can write \[ g_{ij}(x)=\frac{k\phi_{ij}(x)}{x}, \] 
where $\phi_{ij}(x)\phi_{ij}(-x)=1.$ From the results in \cite[Theorem 16]{BFV}, we can conclude that 
\[ g_{ij}(x)=\pm \frac{ke^{\lambda_{ij}x}}{x}, \] where $\lambda_{12}+\lambda_{23}+\lambda_{31}=0.$
When $R$ contains $B_2$ as a subsystem, the argument works well. 
If $R$ contains the subsystem $\{\pm 12,\pm \overline{12},\pm 1, \pm 2 \}$ of type $B_2,$ 
we have the functional equations 
\begin{eqnarray} 
f_{12}(x-y)f_1(x)-f_2(y)f_{12}(x-y)+f_{\overline{12}}(x+y)f_2(y) + f_1(x)f_{\overline{12}}(x+y)=0, \\ 
g_{12}(x-y)g_1(y)-g_2(y)g_{12}(x+y)+g_{\overline{12}}(x+y)g_2(-x) + g_1(x)g_{\overline{12}}(-x+y)=0. 
\end{eqnarray} 
Since $R$ is not of type $A_1$ or $B_2,$ $R$ contains a subsystem of type $A_2.$ 
We may assume that $f_{12},$ $f_{\overline{12}},$ $g_{12}$ and $g_{\overline{12}}$ 
are determined from the subsystems of type $A_2$ in $R$ as follows: 
\[ f_{12}(x)=f_{\overline{12}}(x)=\frac{k}{x}, \; \; g_{12}(x)=\frac{ke^{\lambda_{12}x}}{x}, 
\; \; g_{\overline{12}}(x)= \frac{ke^{\lambda_{\overline{12}}x}}{x}. \] 
Then the functional equations (3) and (4) can be written as 
\begin{eqnarray} 
\left( \frac{1}{x-y}+\frac{1}{x+y} \right)f_1(x) + \left( -\frac{1}{x-y}+\frac{1}{x+y} \right)f_2(y) =0, \\ 
\frac{e^{\lambda_{12}(x-y)}}{x-y}g_1(y)-
\frac{e^{\lambda_{12}(x+y)}}{x+y} g_2(y) + \frac{e^{\lambda_{\overline{12}}(x+y)}}{x+y}g_2(-x)+
\frac{e^{\lambda_{\overline{12}}(-x+y)}}{-x+y} g_1(x) =0. 
\end{eqnarray} 
Hence we get 
\[ f_1(x)=f_2(x)=\frac{k'}{x} \] 
from the equation (5). The equation (6) is written as 
\[ (x+y)e^{\lambda_{12}(x-y)}\frac{\phi_1(y)}{y}-(x-y)e^{\lambda_{12}(x+y)}\frac{\phi_2(y)}{y} \] 
\[ - (x-y)e^{\lambda_{\overline{12}}(x+y)}\frac{\phi_2(-x)}{x}- 
(x+y)e^{\lambda_{\overline{12}}(-x+y)}\frac{\phi_1(x)}{x} \] 
\begin{eqnarray}
= 0.  \;\;\;\;\;\;\;\;\;\;\;\;\;\;\;\;\;\;\;\;\;\;\;\;\;\;\;\;\;\;\;\;\;\;\;\;\;\;\;\;
\;\;\;\;\;\;\;\;\;\;\;\;\;
\end{eqnarray}
We obtain, by taking the limit $y \rightarrow 0,$ 
\[ e^{-\lambda_{\overline{12}}x}\phi_1(x)+e^{\lambda_{\overline{12}}x}\phi_2(-x)
=e^{\lambda_{12}x}(2+x(\phi'_1(0)-\phi'_2(0)-2\lambda_{12})), \] 
and by taking the limit $x \rightarrow 0,$ 
\[ e^{-\lambda_{12}y}\phi_1(y)+e^{\lambda_{12}y}\phi_2(y)=e^{\lambda_{\overline{12}}y}
(2+y(\phi'_1(0)+\phi'_2(0)-2\lambda_{\overline{12}})). \] 
After eliminating $\phi_2(y)$ and $\phi_2(-x)$ from the equation (7), we have 
\[ e^{-(\lambda_{12}+\lambda_{\overline{12}})x}\frac{\phi_1(x)}{x^2}-\frac{1}{x^2}- 
\frac{\phi'_1(0)-\lambda_{12}-\lambda_{\overline{12}}}{x} \] 
\[ =e^{-(\lambda_{12}+\lambda_{\overline{12}})y}\frac{\phi_1(y)}{y^2}-\frac{1}{y^2}- 
\frac{\phi'_1(0)-\lambda_{12}-\lambda_{\overline{12}}}{y}. \] 
This means that the both sides must be a constant $C.$ 
Hence, we have 
\[ \phi_1(x)=e^{(\lambda_{12}+\lambda_{\overline{12}})x}
(1+(\phi'_1(0)-\lambda_{12}-\lambda_{\overline{12}})x+Cx^2) . \] From the condition 
$\phi_1(x)\phi_1(-x)=1,$ we get 
\[ \phi'_1(0)=\lambda_{12}+\lambda_{\overline{12}}, \; \; C=0. \] 
Therefore we conclude that 
\[ g_1(x)=\pm \frac{k'e^{\lambda_1 x}}{x}, \; g_2(x)= \pm 
\frac{k'e^{\lambda_2 x}}{x}, \] where $\lambda_1=\lambda_{12}+\lambda_{\overline{12}},$
$\lambda_2=-\lambda_{12}+\lambda_{\overline{12}}.$ 

When $R_+ = \{ \alpha_1,\alpha_1+\alpha_2, 2\alpha_1+3\alpha_2, \alpha_1+2\alpha_2, \alpha_1+3\alpha_2, \alpha_2 \} $ 
is of type $G_2,$ we have the quadratic relation in the algebra $\B_W$ as follows: 
\[ [\alpha_1][\alpha_1+\alpha_2]+[\alpha_1+\alpha_2][2\alpha_1+3\alpha_2] +[2\alpha_1+3\alpha_2][\alpha_1+2\alpha_2] \] 
\[ + [\alpha_1+2\alpha_2][\alpha_1+3\alpha_2] + [\alpha_1 + 3\alpha_2][\alpha_2] =[\alpha_2][\alpha_1] . \] 
This equation shows that the constants $(\lambda_{\gamma})_{\gamma \in R_+}$ 
are subject to the following constraints 
\[ \lambda_{\alpha_1+\alpha_2}=3\lambda_{\alpha_1}+\lambda_{\alpha_2},
\lambda_{2\alpha_1+3\alpha_2}= 2\lambda_{\alpha_1}+\lambda_{\alpha_2},
\lambda_{\alpha_1+2\alpha_2}=3\lambda_{\alpha_1}+2\lambda_{\alpha_2},
\lambda_{\alpha_1+3\alpha_2}= \lambda_{\alpha_1}+\lambda_{\alpha_2} . \] 
This means that $\lambda_{\gamma}=\lambda(\gamma^{\vee})$ for some $\lambda\in V^*.$ 

Consider the multiplication operators 
\[ {\bf e} = e^{\sum_{i=1}^r \lambda_{\alpha_i}\pi_i(\xi)} \]  
and 
\[ {\bf e}_+=(\prod_{\beta \in R_+}\beta)e^{\sum_{i=1}^r \lambda_{\alpha_i}\pi_i(\xi)}, \]  
where $\pi_i$ is the fundamental dominant weight corresponding to $\alpha_i.$ 
For the operator 
$\DD_{\alpha}=k_{\alpha}(1-e^{\lambda_{\alpha}(\alpha,\xi)}s_{\alpha})/(\alpha,\xi),$ 
we have
\[ \DD_{\alpha}= k_{\alpha} {\bf e} \circ \partial_{\alpha} \circ {\bf e}^{-1}. \] 
For the operator 
$\DD_{\alpha}=k_{\alpha}(1+e^{\lambda_{\alpha}(\alpha,\xi)}s_{\alpha})/(\alpha,\xi),$ 
we have 
\[ \DD_{\alpha}= k_{\alpha} {\bf e}_+ \circ \partial_{\alpha} \circ {\bf e}_+^{-1}. \] 
Namely, 
$\DD_{\alpha}$ is conjugate to $\partial_{\alpha}$ up to a constant $k_{\alpha}.$ 
Hence the operators $\DD_{\alpha}$ 
give rise to a representation of $\B_W$ from Lemma 1.1. 
\begin{prop}
If $R$ is of type $B_2,$ then $\B_W$ is a $64$-dimensional algebra defined 
by the following relations: \\ 
{\rm (i)} $\; [12]^2=\overline{[12]}^2=[1]^2=[2]^2=0,$ \\ 
{\rm (ii)} $\; [12]\overline{[12]}=\overline{[12]}[12],$ $[1][2]=[2][1],$ \\ 
{\rm (iii)} $\; [12][1]-[2][12]+[1]\overline{[12]}+\overline{[12]}[2]=0,$ $[1][12]-[12][2]+
\overline{[12]}[1]+[2]\overline{[12]}=0,$ \\
{\rm (iv)} $\; [12][1]\overline{[12]}[1]+\overline{[12]}[1][12][1]+[1][12][1]\overline{[12]}+ 
[1]\overline{[12]}[1][12]=0,$ \\ 
{\rm (v)} $\; [1][12][1][12]=[12][1][12][1].$
\end{prop} 
The relations above were considered in \cite{KM1} and \cite{MS}. The algebra defined by 
these relations is a finite-dimensional algebra with the Hilbert polynomial 
$(1+t)^4(1+t^2)^2.$ Milinski and Schneider \cite{MS} and Bazlov \cite{Ba} have shown that 
these relations are also satisfied in the algebra 
$\B_W.$ They also checked that the algebra $\B_W$ has dimension 64. 
Hence, the relations above exhaust the independent defining relations 
for the algebra $\B_W$ in $B_2$ case. 
\begin{prop}
Let $R$ be of type $B_2.$ \\  
{\rm (i)} The functions $f_{\alpha}$ must be as follows:  
\[ f_1(x)=f_2(x)=\frac{A}{{\rm sn}(ax,k)}, \] 
\[ f_{12}(x)=f_{\overline{12}}(x)=\frac{B}{{\rm sn}(\varepsilon a x,\tilde{k})}, \] 
where $A,$ $B,$ $a,$ $k$ are arbitrary constants, and 
$\tilde{k}=(1-k)/(1+k),$ $\varepsilon =(1+k)/\sqrt{-1}.$ \\ 
{\rm (ii)} If one assumes the $W$-invariance $w \circ \DD_{\alpha} \circ w^{-1} = \DD_{w(\alpha)},$ 
$w\in W,$ then  
\[ g_{\alpha}(x)=\pm f_{\alpha}(x) , \] 
where the choice of the signature is independent of $\alpha.$ \\ 
{\rm (iii)} If the functions $f_{\alpha}(z)$ are chosen as in {\rm (i)} and 
$g_{\alpha}(x)= \pm e^{\lambda(\alpha^{\vee})x}f_{\alpha}(x),$ 
$\lambda \in V^*,$ then the operators $\DD_{\alpha}$ give a representation of $\B_W.$ 
\end{prop} 
{\it Proof.} (i) This follows from the 4-term quadratic equations and 
\cite[Theorem 6]{BFV}. 
The relation
\[ [12][1]-[2][12]+\overline{[12]}[2]+[1]\overline{[12]}=0 \] 
implies 
\[ f_{12}(x-y)f_1(x)-f_2(y)f_{12}(x-y)+f_{\overline{12}}(x+y)f_2(y)+f_1(x)f_{\overline{12}}(x+y)=0. \] From the equations 
\[ (f_{12}(x-y)+f_{\overline{12}}(-x+y))g_1(x)=0,
\;(f_1(y)-f_2(y))g_{12}(x-y)=0, \]
we have $f_{12}=f_{\overline{12}}$ and
$f_1=f_2.$ Hence the functions $f_{12}=f_{\overline{12}}$ and $f_1=f_2$ are the solutions 
found in \cite{BFV} 
in the invariant case, i.e.,  
\[ f_1(x)=f_2(x)=\frac{A}{{\rm sn}(ax,k)}, \] 
\[ f_{12}(x)=f_{\overline{12}}(x)=\frac{B}{{\rm sn}(\varepsilon a x,\tilde{k})}. \] 
(ii) The $W$-invariance shows $g_{12}(z)=g_{\overline{12}}(z)$ and $g_1(z)=g_2(z).$ 
Moreover, the functions $g_{\alpha}(z)$ must be odd functions. 
On the other hand, we may set $g_{\alpha}(z)=f_{\alpha}(z)\phi_{\alpha}(z)$ with 
$\phi_{\alpha}(z)\phi_{\alpha}(-z)=1$ from Lemma 1.2. Since both of $f_{\alpha}$ and 
$g_{\alpha}$ are odd functions, $\phi_{\alpha}$ must be even function. 
Hence, we have $\phi_{\alpha}(z)=\pm 1.$ \\ 
(iii) In this case, we can check that the operators $\DD_{\alpha}$ satisfy 
all the relations listed in Proposition 1.2 by direct computation. 

\section{Representation of quadratic algebra} 
\begin{dfn} {\rm For a given family of functions $\varphi_{ij}(z)=-\varphi_{ji}(-z),$ 
$\psi_{ij}(z)=\psi_{ji}(z)\in \M_0,$ 
$1\leq i,j \leq n,$ $i\not=j,$ 
the algebra $\tilde{\E}_n(\varphi_{ij},\psi_{ij})$  
is a $\C$-algebra generated by the symbols $\la ij\ra$ and functions $f(\xi)$ in $\M$ subject to the relations: \\ 
(i) $\la ij \ra^2=\psi_{ij}(x_i-x_j),$ \\ 
(ii) $\la ij \ra \la kl \ra = \la kl \ra \la ij\ra$ for $\{ i,j \} \cap \{ k,l \} = \emptyset,$ \\ 
(iii) $\la ij \ra \la jk \ra + \la jk \ra \la ki\ra + \la ki \ra \la ij \ra=0,$ \\ 
(iv) $(\la ij \ra -\varphi_{ij}(x_i-x_j)) f(\xi)=
f(s_{ij}\xi)(\la ij \ra-\varphi_{ij}(x_i-x_j)).$ } 
\end{dfn} 
\begin{rem}
{\rm The algebra $\M^{S_n}$ of $S_n$-invariant functions is contained in the 
center of $\tilde{\E}_n(\varphi_{ij},\psi_{ij}).$ Hence 
$\tilde{\E}_n(\varphi_{ij},\psi_{ij})$ has a structure 
of the $\M^{S_n}$-algebra.} 
\end{rem}
In this section we consider when the quadratic algebra $\tilde{\E}_n(\varphi_{ij},\psi_{ij})$ 
has a generalization of the Calogero-Moser representation. 
For $\lambda \in \C \setminus \Z+\Z \tau,$ 
define the function $\sigma_{\lambda}(z)=\sigma_{\lambda}(z|\tau)$ by the formula 
\[ \sigma_{\lambda}(z)=
\frac{\vartheta_1(z-\lambda)\vartheta_1'(0)}{\vartheta_1(z)\vartheta_1(-\lambda)}, \] 
where $\vartheta_1(z)$ is Jacobi's theta function 
\[ \vartheta_1(z)= -\sum_{n=-\infty}^{+\infty}\exp \left( 2\pi \sqrt{-1} \left( 
(z+\frac{1}{2})(n+\frac{1}{2}) + \frac{\tau}{2} (n+\frac{1}{2})^2 \right) \right) . \] 
\begin{prop}
The algebra $\tilde{\E}_n(\varphi_{ij},\psi_{ij})$ has the generalized Calogero-Moser representation 
if and only if $\varphi_{ij}(z)=a/z$ and 
the functions $\psi_{ij}$ have one of the following forms: \\ 
{\rm (i)}  
\[ \psi_{ij}(z)= \frac{A}{z^2}-K(\wp(bz)-\wp(\lambda_i-\lambda_j)), \] 
{\rm (ii)} 
\[ \psi_{ij}(z)= \frac{A}{z^2}-K\frac{\sin^2(b(z-\lambda_i+\lambda_j))}{\sin^2(bz)\sin^2(b(\lambda_i-\lambda_j))}, \] 
{\rm (iii)} 
\[ \psi_{ij}(z)= \frac{A-K}{z^2}+\frac{K}{(\lambda_i-\lambda_j)^2}. \] 
Here, $A=a^2,$ $K$ and $b$ are parameters. 
\end{prop} 
{\it Proof.} If the generalized Calogero-Moser representation 
\[ \la ij \ra \mapsto \DD_{ij}= f_{ij}(x_i-x_j)+g_{ij}(x_i-x_j)s_{ij} \] 
is well-defined for the algebra $\tilde{\E}_n(\varphi_{ij},\psi_{ij}),$ 
then $\varphi_{ij}(z)=f_{ij}(z)$ must be a rational function $a/z$ as we have seen in the proof of 
Proposition 1.1. The functions $g_{ij}$ are also determined from \cite[Theorem 16]{BFV}. 
Hence the operator $\DD_{ij}$ 
must be one of the following: \\ 
(i) 
\[ \DD_{ij}= \frac{a}{x_i-x_j} + k \sigma_{\lambda_i-\lambda_j}(b(x_i-x_j))
e^{(\alpha_i-\alpha_j)(x_i-x_j)}s_{ij} , \] 
(ii) 
\[ \DD_{ij}= \frac{a}{x_i-x_j} + k \frac{\sin(b(x_i-x_j-\lambda_i+\lambda_j))}{\sin(b(x_i-x_j))\sin(b(\lambda_i-\lambda_j))}e^{(\alpha_i-\alpha_j)(x_i-x_j)}s_{ij}, \] 
(iii) 
\[ \DD_{ij}= \frac{a}{x_i-x_j} + k \left(\frac{1}{x_i-x_j}-\frac{1}{\lambda_i-\lambda_j} \right)e^{(\alpha_i-\alpha_j)(x_i-x_j)}s_{ij}. \] 
In case (i), we have 
\[ \psi_{ij}(x_i-x_j)= \DD_{ij}^2= \frac{A}{(x_i-x_j)^2}-K(\wp(b(x_i-x_j))-\wp(\lambda_i -\lambda_j)) \]
with $A=a^2,$ $K=k^2.$ We also have the desired result in cases (i) and (ii) in a similar way. 
\begin{rem}{\rm 
The trigonometric solution (ii) is obtained from the elliptic solution (i)
by taking the limit $\tau \rightarrow +\infty \sqrt{-1}$ and replacing $\lambda_i$ by $b\lambda_i.$ 
The rational solution (iii) is obtained from the trigonometric solution by taking the limit 
$b\rightarrow 0$ after replacing $K$ by $Kb^2.$}
\end{rem}

Under the assumption of Proposition 2.1, the functions $\varphi_{ij}(z)$ are determined 
to be the rational function $a/z.$ In the rest of this paper, we denote just by 
$\tilde{\E}_n(\psi_{ij})$ the quadratic algebra $\tilde{\E}_n(\varphi_{ij},\psi_{ij})$ 
with $\varphi_{ij}(z)=a/z.$ If we introduce a new set of generators 
$[ij]=\la ij \ra -a/(x_i-x_j),$ 
then the algebra $\tilde{\E}_n(\psi_{ij})$ is defined by the following relations: \\ 
${\rm (i)}'$ $[ij]^2=\psi_{ij}(x_i-x_j)-A/(x_i-x_j)^2,$ \\ 
${\rm (ii)}'$ $[ij][kl]=[kl][ij]$ for $\{ i,j \} \cap \{ k,l \} = \emptyset,$ \\ 
${\rm (iii)}'$ $[ij][jk]+[jk][ki]+[ki][ij]=0,$ \\ 
${\rm (iv)}'$ $[ij]f(\xi)=f(s_{ij}\xi)[ij].$ 

\section{Subalgebra generated by Dunkl elements} 
In this section, the functions $\psi_{ij}$ are assumed to be chosen as in 
Proposition 2.1 (i) with $K=b=1$ for simplicity. 

We define the Dunkl elements $\theta_i$ in the algebra $\tilde{\E}_n(\psi_{ij})$ 
by the formula 
\[ \theta_i =\sum_{j\not= i} [ij] . \] 
We can easily see the following from the defining quadratic relations for 
$\tilde{\E}_n(\psi_{ij}).$ 
\begin{prop}
The Dunkl elements $\theta_1,\ldots,\theta_n$ commute pairwise. 
\end{prop}

In the rest of this section, we discuss the structure of the commutative subalgebra 
generated by the Dunkl elements $\theta_1,\ldots, \theta_n$ over $\M^{S_n}$ in the algebra 
$\tilde{\E}_n(\psi_{ij}).$ We use an abbreviation $x_{ij}:=x_i-x_j,$ 
$\lambda_{ij}:=\lambda_i-\lambda_j$ in the following. 
\begin{lem}{\rm (\cite[Lemma 7.3]{FK})} 
For distinct $i_1,\ldots,i_k,$ one has the following relation in the algebra 
$\tilde{\E}_n(\psi_{ij})$ for $k\geq 3.$ 
\begin{equation} 
\sum_{a=1}^k [i_a \; i_{a+1}][i_a \; i_{a+2}] \cdots [i_a \; i_k] \cdot [i_a \; i_1] [i_a \; i_2] 
\cdots [i_a \; i_{a-1}] = 0. 
\end{equation}  
\end{lem}
{\it Proof.} The proof is done by induction on $k.$ 
For $k=3,$ the relation (8) is just the 3-term relation 
\[ [i_1\; i_2][i_2 \; i_3]+[i_2\; i_3][i_3\; i_1]+[i_3\; i_1][i_1\; i_2]=0 . \] 
Let $Q_k(i_1,\ldots,i_k)$ denote 
the left-hand side of the above relation. 
By using the 3-term relation 
\[ [i_a\; i_{k-1}][i_a\; i_k]=[i_{k-1}\; i_k][i_a\; i_{k-1}]-[i_a\; i_k][i_{k-1}\; i_k] , \] 
we get 
\begin{eqnarray*}
\lefteqn{\sum_{a=1}^k [i_a \; i_{a+1}][i_a \; i_{a+2}] \cdots [i_a \; i_k] \cdot [i_a \; i_1] [i_a \; i_2] 
\cdots [i_a \; i_{a-1}] } \\  
&=& \sum_{a=1}^{k-2} [i_a \; i_{a+1}] \cdots [i_a \; i_{k-2}] \cdot \Big( 
[i_{k-1}\; i_k][i_a\; i_{k-1}]-[i_a\; i_k][i_{k-1}\; i_k] \Big)
\cdot [i_a \; i_1] 
\cdots [i_a \; i_{a-1}] \\ 
& & {} + [i_{k-1}\; i_k][i_{k-1}\; i_1]\cdots [i_{k-1}\; i_{k-2}] +
[i_k\; i_1][i_k\; i_2] \cdots [i_k\; i_{k-1}] \\ 
&=& [i_{k-1}\; i_k] Q_{k-1}(i_1,\ldots , i_{k-1})-Q_{k-1}(i_1,\ldots,i_{k-2},i_k)[i_{k-1}\; i_k] =0.  
\end{eqnarray*}

\begin{lem}
For distinct $i_1,\ldots,i_k,m,$ one has the following relation in the algebra 
$\tilde{\E}_n(\psi_{ij})$ for $k\geq 2.$ 
\[ 
(-1)^{k+1}\sum_{a=1}^k [i_a \; m][i_{a+1} \; m] \cdots [i_k \; m] \cdot [i_1 \; m] [i_2 \; m] 
\cdots [i_{a-1} \; m] [i_a \; m]  \] 
\begin{equation} 
= \sum_{a=1}^k \wp(\lambda_{i_am}) [i_a \; i_{a+1}][i_a \; i_{a+2}] 
\cdots [i_a \; i_k ] \cdot [i_a \; i_1] [i_a \; i_2] \cdots [i_a \; i_{a-1}] 
\end{equation} 
\end{lem}
{\it Proof.} The proof is done by induction on $k.$ 
For $k=2,$ we have 
\begin{eqnarray*}
\lefteqn{[i_1\; m][i_2\; m][i_1\; m]+[i_2\; m][i_1\; m][i_2\; m]} \\ 
&=& \Big( [i_2\; m][i_1\; i_2]-[i_1\; i_2][i_1\; m] \Big) [i_1\; m] + 
[i_2\; m] \Big( [i_2\; m][i_1\; i_2]-[i_1\; i_2][i_1\; m] \Big) \\  
&=& -[i_1\; i_2](\psi_{i_1 \; m}(x_{i_1 \; m})-Ax_{i_1 \; m}^{-2})+ 
(\psi_{i_2\; m}(x_{i_2\; m})-Ax_{i_2\; m}^{-2})[i_1\; i_2] \\ 
&=& \Big( \psi_{i_2\; m}(x_{i_2\; m}) - \psi_{i_1 \; m}(x_{i_2\; m})\Big) [i_1\; i_2] \\ 
&=& \Big( \wp(\lambda_{i_2\; m}) - \wp(\lambda_{i_1\; m}) \Big) [i_1\; i_2] . 
\end{eqnarray*} 
Let $P_k(i_1,\ldots,i_k;m)$ denote 
the left-hand side of the relation (9). 
Here we show only the relation 
\[ P_k(1,2,\ldots,k;m) = \sum_{a=1}^k \wp(\lambda_{am}) [a \; a+1][a \; a+2] 
\cdots [a \; k ] \cdot [a \; 1] [a \; 2] \cdots [a \; a-1], \] 
since the general relations can be proved in similar manner. 
By using the quadratic relation 
$[i_{k-1} \; m][i_k\; m] = [i_k\; m][i_{k-1}\; i_k] -[i_{k-1}\; i_k][i_{k-1}\; m]$ and 
the assumption of the induction, 
we obtain 
\begin{eqnarray*} 
\lefteqn{P_k(1,\ldots,k;m)}\\
&=& [k-1\; k] \cdot P_{k-1}(1,\ldots,k-2,k-1;m)-
P_{k-1}(1,\ldots,k-2,k;m)\cdot [k-1\; k] \\ 
&=& [k-1\; k] \cdot \sum_{a=1}^{k-1} \wp(\lambda_{am}) [a \; a+1]
\cdots [a \; k-1 ] \cdot [a \; 1]  \cdots [a \; a-1] \\ 
& & {} - \sum_{a=1}^{k-2} \wp(\lambda_{am}) [a \; a+1]
\cdots [a \; k-2 ] [a \; k] \cdot [a \; 1]  \cdots [a \; a-1] \cdot [k-1 \; k] \\ 
& & {} - \wp(\lambda_{km}) [k \; 1][k\; 2] \cdots [k \; k-2] [k-1 \; k] \\ 
&=& \sum_{a=1}^{k-2} \wp(\lambda_{am}) [a \; a+1]
\cdots [a \; k-2 ] \Big( [k-1 \; k][a \; k-1] - [a \; k][k-1 \; k] \Big) [a \; 1]  \cdots [a \; a-1] \\ 
& & {} + \wp(\lambda_{k-1 \; m})[k-1 \; k][k-1\; 1] \cdots [k-1 \; k-2] \\ 
& & {} + \wp(\lambda_{km}) [k \; 1][k\; 2] \cdots [k \; k-2] [k\; k-1] \\ 
&=& \sum_{a=1}^{k-2} \wp(\lambda_{am}) [a \; a+1]
\cdots [a \; k-2 ] \Big( [a \; k-1][a \; k] \Big) [a \; 1]  \cdots [a \; a-1] \\ 
& & {} + \wp(\lambda_{k-1 \; m})[k-1 \; k][k-1\; 1] \cdots [k-1 \; k-2] \\ 
& & {} + \wp(\lambda_{km}) [k \; 1][k\; 2] \cdots [k \; k-2] [k\; k-1] \\ 
&=& \sum_{a=1}^k \wp(\lambda_{am}) [a \; a+1][a \; a+2] 
\cdots [a \; k ] \cdot [a \; 1] [a \; 2] \cdots [a \; a-1]. 
\end{eqnarray*}  
\\ 
{\bf Example.} $(k=4)$ 
\begin{eqnarray*}
\lefteqn{[1m][2m][3m][4m][1m]+ [2m][3m][4m][1m][2m]}\\ 
\lefteqn{+[3m][4m][1m][2m][3m]+ [4m][1m][2m][3m][4m]} \\
&=& -[1m][2m][34][3m][1m]+[1m][2m][4m][34][1m] \\ 
& & {} -[2m][34][3m][1m][2m]+[2m][4m][34][1m][2m] \\ 
& & {} -[34][3m][1m][2m][3m]+[4m][34][1m][2m][3m] \\ 
& & {} -[4m][1m][2m][34][3m]+[4m][1m][2m][4m][34] \\ 
&=& - [34] \Big( [1m][2m][3m][1m]+[2m][3m][1m][2m]+[3m][1m][2m][3m] \Big) \\ 
& & {} + \Big( [1m][2m][4m][1m]+[2m][4m][1m][2m]+[4m][1m][2m][4m] \Big) [34] \\  
&=& -[34] \Big( \wp(\lambda_{1m}) [12][13] + 
\wp(\lambda_{2m}) [23][21] + \wp(\lambda_{3m}) [31][32] \Big) \\ 
& & + \Big( \wp(\lambda_{1m}) [12][14] + 
\wp(\lambda_{2m}) [24][21] + \wp(\lambda_{4m}) [41][42] \Big) [34] \\ 
&=& -\wp(\lambda_{1m}) [12] ([34][13]-[14][34]) - 
\wp(\lambda_{2m})([34][23]-[24][34])[21] \\ 
& & - \wp(\lambda_{3m}) [34][31][32] - \wp(\lambda_{4m})[41][42][43] \\ 
&=& -\wp(\lambda_{1m}) [12] [13][14] - \wp(\lambda_{2m}) [23][24][21] - 
\wp(\lambda_{3m}) [34] [31][32] - \wp(\lambda_{4m}) [41] [42][43]  
\end{eqnarray*}
\begin{rem}
{\rm Lemma 3.2 is a deformed version of \cite[Lemma 7.2]{FK} and \cite[Lemma 5.3]{Po}. 
Though the identity (9) looks similar to the one in \cite[Lemma 5.3]{Po}, 
they are different formulas. In our case, $[ij]^2=\psi_{ij}(x_{ij})-Ax_{ij}^{-2}$ is not central, 
and $[ij]^2\not= \wp(\lambda_{ij}).$}
\end{rem}

For a subset $I\subset \{ 1,\ldots ,n \}$ with $\# I=2k,$ define the function 
$\phi(I)=\phi(x_i|i\in I)$ by the following formula: 
\[ \phi(I):= \sum \prod_{i=1}^k \wp(x_{a_i b_i}) , \] 
where the summension is taken over the choice of pairs $(a_i,b_i),$ $1\leq i \leq k,$ 
such that $I=\{ a_1,\ldots,a_k, b_1,\ldots b_k \},$ $a_1< \cdots < a_k$ and $a_i<b_i.$ 
We also define the deformed elementary symmetric polynomial $E_k(I)=E_k(X_i\, |\, i \in I)$ 
by the recursion relations: 
\[ E_0(I)=1, \; E_k(I\cup \{ j \})=E_k(I)+E_{k-1}(I)X_j+\sum_{i\in I}
\wp(\lambda_{ij})E_{k-2}(I\setminus \{ i \}). \]
\begin{thm}
One has the following formula in the algebra $\tilde{\E}_n(\psi_{ij}):$ 
\begin{equation} 
E_k(\theta_i \; | \; i\in I)= \sum_{l=0}^{[k/2]}\sum_{I_0 \subset I, \# I_0=2l} \phi(I_0) \sum_{(*)} 
[a_1 \; b_1] \cdots [a_{k-2l} \; b_{k-2l}] , 
\end{equation} 
where $(*)$ stands for the conditions that $a_i \in I \setminus I_0$; $b_i \not\in I$; 
$a_1,\ldots ,a_{k-2l}$ are distinct; $b_1\leq \cdots \leq b_{k-2l}.$ 
\end{thm}
\begin{cor} {\rm 
\[ E_k(\theta_1,\ldots,\theta_n)= \left\{ 
\begin{array}{cc}
\sum_{I_0\subset I, \# I_0=k} \phi(I_0) & \textrm{if $k$ is even,} \\ 
0 & \textrm{if $k$ is odd.} 
\end{array}
\right. \] }
\end{cor}
{\it Proof of Theorem 3.1.} Denote by $F_k(I)$ the right-hand side of the fomula (10). 
For $I\subset \{ 1, \ldots, n\}$ and $j\not\in I,$ we will show the recursive formula 
\[ F_k(I\cup \{ j \})=F_k(I)+\theta_jF_{k-1}(I)+\sum_{i\in I} 
\wp(\lambda_{ij})F_{k-2}(I\setminus \{ i \}). \] 
Let $J=\{ j_1=j,\ldots ,j_d \}$ be the set $\{ 1,\ldots , n \} \setminus I.$ For 
$L=\{ l_1, \ldots, l_m \} \subset \{ 1,\ldots , n \}$ and $r \not\in L,$ we define 
\[ \lla L \; | \; r \rra := 
\sum_{w \in S_m} [l_{w(1)} \; r][l_{w(2)} \; r] \cdots [l_{w(m)} \; r] . \] 
In order to show the formula above, we use 
the following decompositions which are 
similar to those used in \cite{Po}. In the following, the symbol 
$I_1\cdots I_d \subset_m I$ means that $I_1,\ldots,I_d \subset I$ are disjoint and 
$\# I_1 + \cdots + \# I_d=m.$ Here, some of $I_1, \ldots, I_d$ may be empty. 
Let us consider the decompositions: 
\begin{eqnarray*} 
F_k(I) &=& \sum_{l=0}^{[k/2]}\sum_{I_0\subset_{2l} I}\phi(I_0)
\sum_{I_1\ldots I_d \subset_{k-2l} I\setminus I_0} \lla I_1 \; | \; j_1 \rra 
\lla I_2 \; | \; j_2 \rra \cdots \lla I_d \; | \; j_d \rra \\ 
&=& A_1 + A_2, \\ 
F_k(I\cup \{ j \}) &=& \sum_{l=0}^{[k/2]}\sum_{I'_0\subset_{2l} I\cup \{ j \} }\!\!\!\! 
\phi(I'_0)\!\!\!\!\! \sum_{I'_2\ldots I'_d \subset_{k-2l} (I\cup \{ j \} ) \setminus I'_0} \!\!\!\!\!\!\!\!\!\!\!\!\! 
\lla I'_2 \; | \; j_2 \rra \lla I'_3 \; | \; j_3 \rra \cdots \lla I'_d \; | \; j_d \rra \\ 
&=& B_1 + B_2 + B_3, \\ 
\theta_jF_{k-1}(I) &=& \sum_{s\not=j}[js]\sum_{l=0}^{[(k-1)/2]}\sum_{I''_0\subset_{2l} I}\!\!\! \phi(I''_0)
\!\!\! \sum_{I''_1\ldots I''_d \subset_{k-1-2l} I\setminus I''_0} \!\!\!\!\!\! \lla I''_1 \; | \; j_1 \rra 
\cdots \lla I''_d \; | \; j_d \rra \\ 
&=& C_1+C_2+C_3+C_4, 
\end{eqnarray*}  
where $A_1$ is the sum of terms with $I_1=\emptyset$; $A_2$ is the sum of terms with 
$I_1\not= \emptyset$; $B_1$ is the sum of the terms with $j\not\in I_0 \cup I'_2 \cup \cdots 
\cup I'_d$; $B_2$ is the sum of terms with $j\in I'_2 \cup \cdots \cup I'_d$; 
$B_3$ is the sum of terms with $j\in I''_0$; $C_1$ is the sum of terms with 
$s\in I\setminus (I''_0\cup I''_1 \cup \cdots \cup I''_d)$; $C_2$ is the sum of terms with 
$s\in I''_2 \cup \cdots \cup I''_d \cup J$; $C_3$ is the sum 
of terms with $s\in I''_0$; $C_4$ is the sum of terms with $s\in I''_1.$ Then we can see that 
$A_1=B_1,$ $A_2+C_1=0$ and $B_2=C_2$ by the same argument in \cite{Po}. 

Note that the formula in Lemma 3.2 holds only for $k\geq 2.$ For any subset 
$K=\{ k_1,\ldots,k_m \}$ with $j\not\in K,$ 
we have 
\begin{eqnarray*} 
\lefteqn{\sum_{s\in K}\left( [js] \lla K \; | \; j \rra + 
\sum_{L\subset K\setminus \{ s \}} \wp(\lambda_{js})\lla L \; | \; s \rra \lla K\setminus L \setminus \{ s \} 
\; | \; j \rra   \right)} \\ 
&=& \sum_{s\in K} \left( [js][sj]\lla K\setminus \{ s \} \; | \; j \rra + 
[js] \sum_{w\in S_m, k_{w(1)}\not=s} [ k_{w(1)}\; j] \cdots [k_{w(k)} \; j ] \right. \\ 
& & \left. + \wp(\lambda_{js})\lla K \setminus \{ s \} 
\; | \; j \rra  + \sum_{L\subset K\setminus \{ s \},L\not=\emptyset} 
\wp(\lambda_{js})\lla L \; | \; s \rra \lla K\setminus L \setminus \{ s \} 
\; | \; j \rra   \right) \\
&=& \sum_{s\in K}\wp(x_{js}) \lla K \setminus \{ s \} 
\; | \; j \rra  
\end{eqnarray*}
from Lemma 3.2 and $[ij]^2=\psi_{ij}(x_{ij})-Ax_{ij}^{-2}=
-(\wp(x_{ij})-\wp(\lambda_{ij})).$ 
This shows 
\begin{eqnarray*} 
\lefteqn{C_3+C_4+\sum_{i\in I}\wp(\lambda_{ij})F_{k-2}(I\setminus \{ i \})} \\ 
& = & \sum_{l=1}^{[(k-1)/2]}\!\!\! \sum_{I''_0\subset_{2l} I} \sum_{s\in I''_0} 
\phi(\{ j \}\cup I''_0 \setminus \{ s \})\cdot [js] \!\!
\sum_{I''_1\cdots I''_d \subset_{k-1-2l}I\setminus I_0''}\!\!\!\!\!\! \lla I''_1 \; | \; j_1 \rra 
\cdots \lla I''_d \; | \; j_d \rra \\ 
& + & \!\!\!\!\!\! \sum_{l=0}^{[(k-1)/2]}\!\!\!\!\! \sum_{I''_0\subset_{2l} I}\!\! \phi(I''_0) \!\!\!
\sum_{I''_1\cdots I''_d \subset_{k-1-2l}I\setminus I''_0}\sum_{s\in I''_1} 
\wp(x_{js}) \lla I''_1 \setminus \{ s \}  | \; j_1 \rra 
\cdots \lla I''_d \; | \; j_d \rra \\ 
& = & -\sum_{l=1}^{[(k-1)/2]}\sum_{I''_0\subset_{2l} I,j\in I''_0}\phi(I_0'')
\sum_{I''_1\cdots I''_d \subset_{k-2l}I\setminus I''_0, I''_1\not=\emptyset}
\lla I''_1 \; | \; j_1 \rra \cdots \lla I''_d \; | \; j_d \rra \\ 
& + & \sum_{l=0}^{[(k-1)/2]}\sum_{I''_0\subset_{2l+2} I,j\in I''_0}\phi(I_0'')
\sum_{I''_1\cdots I''_d \subset_{k-2l-2}I\setminus I''_0}
\lla I''_1 \; | \; j_1 \rra \cdots \lla I''_d \; | \; j_d \rra \\ 
& = & \sum_{l=1}^{[k/2]}\sum_{I''_0\subset_{2l} I,j\in I''_0}\phi(I_0'')
\sum_{I''_2\cdots I''_d \subset_{k-2l}I\setminus I''_0}
\lla I''_2 \; | \; j_2 \rra \cdots \lla I''_d \; | \; j_d \rra \\ 
& = & B_3. 
\end{eqnarray*} 
{\bf Example.} 
One has the following formula for $E_3(\theta_1,\theta_2,\theta_3)$ in 
$\tilde{\E}_5(\psi_{ij}):$
\[ \theta_1\theta_2\theta_3+\wp(\lambda_{23})\theta_1+\wp(\lambda_{13})\theta_2 + 
\wp(\lambda_{12})\theta_3= \] 
\[ \sum_{(**)}[a_1\; b_1][a_2\; b_2][a_3\; b_3]+\psi_{12}(x_{12})([34]+[35])+\psi_{13}(x_{13})([24]+[25])+
\psi_{23}(x_{23})([14]+[15]), \] 
where $(**)$ stands for the condition that $\{ a_1,a_2,a_3 \} =\{ 1,2,3 \};$ 
$b_1,b_2,b_3\in \{ 4,5 \}$ and $b_1\leq b_2 \leq b_3.$ 

\section{Degenerations} 
Some variants of the cohomology ring of the flag variety 
\[ Fl_n= SL_n(\C)/\textrm{(upper triangular matrices)} \] 
have the model 
as the commutative subalgebra in deformations of the quadratic algebra $\E_n.$ 
We see how the deformations of $\E_n$ used for the constructions of the 
model of the cohomology rings can be recovered as degenerations of our algebra 
$\tilde{\E}_n(\psi_{ij}).$ 

Let $T\subset SL_n(\C)$ be the torus consisting of the diagonal matrices. 
We identify the polynomial ring $R=\Z[x_1,\ldots,x_n]$ with the $T$-equivariant 
cohomology ring $H_T(\textrm{pt.}).$ 
The authors introduced the extended quadratic algebra $\tilde{\E}_n \langle R \rangle [t]$ 
to construct a model of the $T$-equivariant cohomology ring $H_T(Fl_n)$ in \cite{KM2}. 
In case $\psi_{ij}(z)=0$ for any distinct $i$ and $j,$ the algebra $\tilde{\E}_n(\psi_{ij}=0)$ 
is defined by the relations $[ij]^2=0,$ $[ij][kl]=[kl][ij]$ for $\{ i,j \} \cap \{ k,l \} = \emptyset,$ 
$[ij][jk]+[jk][ki]+[ki][ij]=0$ 
and $[ij]x_i=x_j[ij].$ Since these relations are same as the 
defining relations for the algebra $\tilde{\E}_n \langle R \rangle [t]|_{t=0}$ introduced in \cite{KM2}, 
the $\C$-subalgebra of 
$\tilde{\E}_n(\psi_{ij}=0)$ generated by $[ij]$'s and $x_1,\ldots,x_n$ 
is isomorphic to $\tilde{\E}_n \langle R \rangle [t]|_{t=0}.$ 
The subsequent result shows that the elements 
\[ \theta'_i:=x_i+\theta_i=x_i+\sum_{j\not=i}[ij], \; \; \; i=1,\ldots,n, \] 
generate a commutative $R$-subalgebra of $\tilde{\E}_n(\psi_{ij}=0)\otimes_{R^{S_n}}R$ 
which is isomorphic to the $T$-equivariant cohomology ring 
$H_T(Fl_n).$ 
\begin{prop} {\rm (\cite[Corollary 2.2]{KM2})}
Let $I$ be a subset of $\{ 1,\ldots,n \}.$ 
In the algebra $\tilde{\E}_n(\psi_{ij}=0),$ one has 
\begin{equation} 
e_k(\theta'_i\; | \; i\in I) = \sum_{m=0}^k\sum_{I_0\subset_m I}(\prod_{i\in I_0}x_i)
\sum_{(*)}[a_1\; b_1]\cdots [a_{k-m}\; b_{k-m}] , 
\end{equation} 
where $(*)$ stands for the conditions that $a_i \in I\setminus I_0$; $b_i\not\in I$; 
$a_1,\ldots,a_{k-m}$ are distinct; $b_1\leq \cdots \leq b_{k-m}.$ 
In particular, one has 
\[ e_k(\theta'_1,\ldots,\theta'_n)=e_k(x_1,\ldots,x_n), \; \; \; 1\leq k \leq n. \] 
\end{prop}
{\it Proof.} The idea is similar to the proof of Theorem 3.1. 
Denote by $F'_k(I)$ the right-hand side of (11). For $j\not\in I,$ we will 
show that 
\[ F'_k(I\cup \{ j \})=F'_k(I)+ F'_{k-1}(I)(x_j+\theta_j). \] 
We use the same notation as the one used in the proof of Theorem 3.1. 
Let us consider the decompositions: 
\begin{eqnarray*} 
F'_k(I) &=& \sum_{m=0}^{k}\sum_{I_0\subset_m I}(\prod_{i\in I_0}x_i)
\sum_{I_1\ldots I_d \subset_{k-m} I\setminus I_0} \lla I_1 \; | \; j_1 \rra 
\lla I_2 \; | \; j_2 \rra \cdots \lla I_d \; | \; j_d \rra \\ 
&=& A'_1 + A'_2, \\ 
F'_k(I\cup \{ j \}) &=& \sum_{m=0}^{k}\sum_{I'_0\subset_m I\cup \{ j \} }\! 
(\prod_{i\in I'_0}x_i)\!\!\!\!\! \sum_{I'_2\ldots I'_d \subset_{k-m} (I\cup \{ j \} ) \setminus I'_0} \!\!\!\!\!\!\!\!\!\!\!\!\! 
\lla I'_2 \; | \; j_2 \rra \lla I'_3 \; | \; j_3 \rra \cdots \lla I'_d \; | \; j_d \rra \\ 
&=& B'_1 + B'_2 + B'_3, \\ 
F'_{k-1}(I)\theta_j &=& \sum_{m=0}^{k-1}\sum_{I''_0\subset_m I}\! (\prod_{i\in I''_0}x_i)
\!\!\! \sum_{I''_1\ldots I''_d \subset_{k-1-m} I\setminus I''_0} \!\!\!\!\!\! \lla I''_1 \; | \; j_1 \rra 
\cdots \lla I''_d \; | \; j_d \rra \sum_{s\not=j}[js] \\ 
&=& C'_1+C'_2+C'_3+C'_4, 
\end{eqnarray*} 
where $A'_1$ is the sum of terms with $I_1=\emptyset$; $A'_2$ is the sum of terms with 
$I_1\not= \emptyset$; $B'_1$ is the sum of the terms with $j\not\in I_0 \cup I'_2 \cup \cdots 
\cup I'_d$; $B'_2$ is the sum of terms with $j\in I'_2 \cup \cdots \cup I'_d$; 
$B'_3$ is the sum of terms with $j\in I''_0$; $C'_1$ is the sum of terms with 
$s\in I\setminus (I''_0\cup I''_1 \cup \cdots \cup I''_d)$; $C'_2$ is the sum of terms with 
$s\in I''_2 \cup \cdots \cup I''_d \cup J$; $C'_3$ is the sum 
of terms with $s\in I''_0$; $C'_4$ is the sum of terms with $s\in I''_1.$ 
Moreover, we decompose $F'_{k-1}(I)x_j$ as follows: 
\begin{eqnarray*} 
F'_{k-1}(I)x_j & = &  \sum_{m=0}^{k-1}\sum_{I''_0\subset_m I} (\prod_{i\in I''_0}x_i)
\!\!\! \sum_{I''_1\ldots I''_d \subset_{k-1-m} I\setminus I''_0} \!\!\!\!\!\! \lla I''_1 \; | \; j_1 \rra 
\cdots \lla I''_d \; | \; j_d \rra x_j \\ 
 &=& D'_1+D'_2, 
\end{eqnarray*} 
where $D'_1$ is the sum of terms with $I''_1=\emptyset$ and $D'_2$ is the sum of terms with 
$I''_1 \not= \emptyset.$ 
As before, we can easily see that $A'_1=B'_1,$ $A'_2+C'_1=0$ and $B'_2=C'_2.$ 
It is also clear that $B'_3=D'_1.$ 
Since the relations $[ij]^2=0$ are assumed, the degenerate version of the formula (9), 
which is same as \cite[Lemma 7.2]{FK}, holds in $\tilde{\E}_n(\psi_{ij}=0):$ 
\[ \sum_{a=1}^k [i_a \; m][i_{a+1} \; m] \cdots [i_k \; m] \cdot [i_1 \; m] [i_2 \; m] 
\cdots [i_{a-1} \; m] [i_a \; m] =0, \; \; \textrm{for $k\geq 1.$} \] 
This formula implies $C'_4=0.$ Finally, the following computation completes the proof: 
\begin{eqnarray*} 
\lefteqn{C'_3+D'_2} \\ 
& = & \sum_{m=1}^{k-1}\sum_{I''_0\subset_m I} (\prod_{i\in I''_0}x_i)
\!\!\! \sum_{I''_1\ldots I''_d \subset_{k-1-m} I\setminus I''_0} \!\!\!\!\!\! \lla I''_1 \; | \; j_1 \rra 
\cdots \lla I''_d \; | \; j_d \rra \sum_{s \in I''_0}[js] \\ 
& + & \sum_{m=0}^{k-1}\sum_{I''_0\subset_m I} (\prod_{i\in I''_0}x_i)
\!\!\! \sum_{I''_1\ldots I''_d \subset_{k-1-m} I\setminus I''_0,I''_1\not=\emptyset} \!\!\!\!\!\! \lla I''_1 \; | \; j_1 \rra 
\cdots \lla I''_d \; | \; j_d \rra x_j \\ 
& = & \!\!\! -\sum_{m=1}^{k-1}\sum_{I''_0\subset_m I} (\prod_{i\in I''_0}x_i)
\!\!\!\! \sum_{I''_1\ldots I''_d \subset_{k-1-m} I\setminus I''_0}  
\sum_{s \in I''_0} \lla I''_1 \; | \; j_1 \rra [sj] \lla I_2 \; | \; j_2 \rra 
\cdots \lla I''_d \; | \; j_d \rra \\ 
 & + & \!\!\!\! \sum_{m=0}^{k-1}\sum_{I''_0\subset_m I} (\prod_{i\in I''_0}x_i)
\!\!\!\!\! \sum_{I''_1\ldots I''_d \subset_{k-1-m} I\setminus I''_0,I''_1\not=\emptyset}  
\sum_{s\in I''_1}x_s\lla I''_1\setminus \! \{ s \} | j_1 \rra [sj] \lla I_2  | j_2 \rra
\cdots \lla I''_d  |  j_d \rra \\ 
 & = & 0. 
\end{eqnarray*}

Let us consider another kind of degeneration. Consider the elliptic solution obtained 
in Proposition 2.1 (i). If we put $K=\kappa \delta^2$ and $\lambda_{ij}=\delta \Lambda_{ij},$ 
then we have 
\[ \lim_{\delta \rightarrow 0} \psi_{ij}(x_{ij})=Ax_{ij}^{-2}+\kappa \Lambda_{ij}^{-2}. \]  
In this situation, the functions $[ij]^2=\psi_{ij}(x_{ij})-Ax_{ij}^{-2}$ 
become central parameters $\kappa \Lambda_{ij}^{-2}.$ Then the $\C$-algebra generated by 
the brackets $[ij]$ in 
$\tilde{\E}_n(\psi_{ij}=\kappa \Lambda_{ij}^{-2})$ is isomorphic to 
the multiparameter deformation of $\E_n$ introduced in \cite[Section 15]{FK}, 
which is denoted by $\E_n^p$ in \cite{Po}, 
under the identification 
of the parameters $p_{ij}=p_{ji}=\kappa \Lambda_{ij}^{-2}.$ In this case, the functions 
$\phi(I)$ are constantly zero, so Theorem 3.1 is reduced to the following: 
\begin{prop} {\rm (\cite[Conjecture 15.1]{FK}, \cite[Theorem 3.1]{Po})} 
Assume that the functions $\psi_{ij}$ are chosen as in Proposition 2.1 (iii) with 
$K=\kappa \delta^2,$ $\lambda_{ij}=\delta \Lambda_{ij}.$ 
In the limit $\delta \rightarrow 0,$ one has
\[ E_k(\theta_i , i\in I ; p)= \sum_{(*)'} 
[a_1 \; b_1] \cdots [a_k \; b_k] , \]
where $(*)'$ stands for the conditions that $a_i \in I$; $b_i \not\in I$; 
$a_1,\ldots ,a_k$ are distinct; $b_1\leq \cdots \leq b_k.$  
\end{prop}

Research Institute for Mathematical Sciences \\ 
Kyoto University \\
Sakyo-ku, Kyoto 606-8502, Japan \\
e-mail: {\tt kirillov@kurims.kyoto-u.ac.jp} \\
URL: {\tt http://www.kurims.kyoto-u.ac.jp/\textasciitilde kirillov} 
\bigskip \\ 
Department of Electrical Engineering, \\
Kyoto University, \\ 
Sakyo-ku, Kyoto 606-8501, Japan \\ 
e-mail: {\tt maeno@kuee.kyoto-u.ac.jp}
\end{document}